\documentclass[12pt]{amsart}
\usepackage{amsfonts,amssymb,amscd,amstext,mathrsfs}
\usepackage[utf8]{inputenc}
\usepackage{hyperref}
\usepackage{verbatim}

\usepackage{graphics}
\usepackage{graphicx}

\usepackage{times}
\usepackage{enumerate}
\usepackage[up,bf]{caption}

\input xy
\xyoption{all}

\newtheorem{theorem}{Theorem}[section]
\newtheorem{proposition}[theorem]{Proposition}

\newtheorem{lemma}[theorem]{Lemma}

\theoremstyle{definition}
\newtheorem{definition}[theorem]{Definition}
\newtheorem{remark}[theorem]{Remark}

\newtheorem{problem}[theorem]{Problem}

\numberwithin{equation}{section}
\numberwithin{figure}{section}

\newcommand\Ascr{\mathscr{A}}

\newcommand\Cscr{\mathscr{C}}

\newcommand\C{\mathbb{C}}
\newcommand\D{\overline{\mathbb D}}
\newcommand\CP{\mathbb{CP}}

\renewcommand\D{\mathbb D}

\newcommand\N{\mathbb{N}}

\newcommand\R{\mathbb{R}}

\newcommand\igot{\mathfrak{i}}

\renewcommand\igot{\mathfrak{i}}

\renewcommand\imath{\igot}

\newcommand\hra{\hookrightarrow}

\newcommand\di{\partial}

\newcommand\dist{\mathrm{dist}}

\newcommand\cd{\overline{\D}}

\def\dist{\mathrm{dist}}

\newcommand\nullq{{\mathbf A}}

\numberwithin{equation}{section}

\begin{document}
\title
{The Calabi--Yau problem for minimal surfaces \\ with Cantor ends}
\author{Franc Forstneri{\v c}}

\address{Franc Forstneri\v c, Faculty of Mathematics and Physics, University of Ljubljana, and Institute of Mathematics, Physics, and Mechanics, Jadranska 19, 1000 Ljubljana, Slovenia}
\email{franc.forstneric@fmf.uni-lj.si}

\thanks{Research is supported by the grants P1-0291, J1-3005, and N1-0237 from ARRS, Republic of Slovenia.}


\date{February 16, 2022. This version July 14, 2022.}


\begin{abstract}
We show that every connected compact or bordered Riemann surface  contains a Cantor set whose 
complement admits a complete conformal minimal immersion in $\R^3$ with bounded image. 
The analogous result holds for holomorphic immersions into any complex manifold of dimension at least $2$,
for holomorphic null immersions into $\C^n$ with $n\ge 3$, for holomorphic Legendrian immersions into
an arbitrary complex contact manifold, and for superminimal immersions into any self-dual or
anti-self-dual Einstein four-manifold.
\end{abstract}

\maketitle


%
%
%
%
\section{Introduction}\label{sec:intro} 
Let $M$ be an open Riemann surface. It is classical \cite{Osserman1986,AlarconForstnericLopez2021} that  
an immersion $M\to\R^n$ for $n\ge 3$ which is conformal (angle preserving) and harmonic 
parameterizes a minimal surface in $\R^n$; conversely, every immersed minimal surface in $\R^n$ arises in this way. 

Let $ds^2$ denote the Euclidean metric on $\R^n$. An immersion $f:M\to \R^n$ is said to be {\em complete} if the 
Riemannian metric $f^* ds^2$ on $M$ induces a complete distance function;  equivalently, if the image of any divergent path in $M$ 
by the map $f$ is a path in $\R^n$ with infinite Euclidean length. 

The {\em Calabi--Yau problem for minimal surfaces} (see \cite[p.\ 170]{Calabi1965Conjecture}
and \cite[p.\ 360]{Yau2000AMS}) asks about the existence,
conformal and asymptotic properties of complete immersed or embedded minimal surfaces 
with bounded images in $\R^n$ for $n\ge 3$. Pioneering constructions 
were given by Jorge and Xavier \cite{JorgeXavier1980AM} in 1980 and Nadirashvili \cite{Nadirashvili1996IM} in 1996.
There were substantial developments since then, and a survey 
can be found in \cite[Chapter 7]{AlarconForstnericLopez2021}.

In this paper we construct the first known examples of complete bounded minimal surfaces whose end is a {\em Cantor set}, that is,
a compact, perfect, totally disconnected set.  The following is a special case of our main result, Theorem \ref{th:mainbis}.
See also Theorem \ref{th:main2} for a generalization to a number of other geometries.

%
%
\begin{theorem}\label{th:main}
In every compact connected Riemann surface, $M$, there is a Cantor set $C$ 
whose complement admits a complete conformal minimal immersion 
$M\setminus C \to\R^3$ with bounded image. There also exist a Cantor set $C$ in $M$ 
and a complete conformal minimal embedding $M\setminus C \hra\R^5$ with bounded image. 
\end{theorem}
 
Cantor sets, being of fractal nature, often serve as a challenging test case in geometric problems.
Theorem \ref{th:main} gives an affirmative answer to Problem 7.4.8 (B) in \cite{AlarconForstnericLopez2021}. 
Note that the problem was posed incorrectly: whether there is a Cantor set $C\subset \C$ such that 
$\C\setminus C$ satisfies the conclusion of the theorem. This is impossible since a bounded harmonic map 
$\C\setminus C\to\R^n$ extends harmonically across the puncture at infinity. 
The correct question is answered affirmatively by Theorem \ref{th:main} with $M=\CP^1$. 

%
%
\begin{remark}\label{rem:main}
Theorem \ref{th:main} also holds if $M$ is a bordered Riemann surface of the form
\begin{equation}\label{eq:bordered}
	M=R\setminus \bigcup_i D_i, 
\end{equation}
where $R$ is a compact connected Riemann surface and $\{D_i\}_i$ is a finite or 
countable collection of pairwise disjoint, smoothly bounded closed discs, 
diffeomorphic images of the unit disc $\cd=\{z\in\C:|z|\le 1\}$. (The discs $D_i$ 
may cluster on one another, but $M$ must be an open domain in $R$.)
It was shown by Alarc\'on and the author \cite{AlarconForstneric2021RMI} in 2021 that 
such $M$ admits a bounded complete conformal minimal immersion in $\R^3$ and embedding 
in $\R^5$ extending continuously to $\overline M$ such that the image of the boundary $bM=\bigcup_i bD_i$ 
is a union of pairwise disjoint Jordan curves, the images of $bD_i$. (For finitely many discs $D_i$ this 
was proved beforehand in \cite{AlarconDrinovecForstnericLopez2015PLMS}.) 
Together with our proof of Theorem \ref{th:main} (see Sect.\ \ref{sec:proof1}), this gives 
a Cantor set $C\subset M$ and a bounded complete conformal 
minimal immersion $M\setminus C\to \R^3$ (and an embedding into $\R^5$)
which extends continuously to $bM=\bigcup_i bD_i$ and maps the curves $bD_i$ to pairwise disjoint Jordan curves.
(However, our proof does not give a continuous extension of the map  to the Cantor set $C$.)  
The point is that we can simultaneously increase 
the intrinsic boundary distances at $bM$ and at $C$. 
We also provide a precise control of the location of the image surface in $\R^n$; see Theorem \ref{th:mainbis}. 
\end{remark}

Cantor sets which arise in the proof Theorem \ref{th:main} are small 
modifications of the standard Cantor sets in the plane, and they have almost full measure in a 
surrounding rectangle. One may ask which Cantor sets in compact Riemann surfaces satisfy
the conclusion of Theorem \ref{th:main}. In my opinion,  this question is likely very difficult or even 
impossible to answer. Instead, I propose the following more reasonable problem.

\begin{problem}
Is there a Cantor set $C$ in $\CP^1$ of zero area such that $\CP^1\setminus C$ satisfies the conclusion of Theorem \ref{th:main}?
\end{problem}

The first test case is to find a Cantor set with zero area whose complement $\CP^1\setminus C$ 
admits a nonconstant bounded harmonic function.

%
%
%
%
%
\section{Tools used in the proof}\label{sec:tools}
In this section we recall the prerequisites and tools which will be used in the proof. They are 
described in detail in the monograph \cite{AlarconForstnericLopez2021}, and we provide precise references.

The following quadric complex hypersurface in $\C^n$ for $n\ge 3$ is called the {\em null quadric}:
\begin{equation} \label{eq:nullq}
	\nullq =\bigl\{z=(z_1,z_2,\ldots,z_n)\in\C^n: z_1^2+z_2^2+\cdots + z_n^2=0\bigr\}.
\end{equation}
Let $M$ be an open Riemann surface. Fix a nowhere vanishing holomorphic 1-form $\theta$ on $M$;
such exists by the Oka--Grauert principle (see \cite[Theorem 5.3.1]{Forstneric2017E}).
An immersion $f=(f_1,\ldots,f_n):M\to\R^n$ is conformal and minimal (equivalently, conformal and harmonic) 
if and only if the $(1,0)$-differential $\di f=(\di f_1,\ldots,\di f_n)$ (the $\C$-linear part of the differential $df$)  
is holomorphic and satisfies the nullity condition $\sum_{i=1}^n (\di f_i)^2=0$. Equivalenty, the map 
\[
	h=2\, \di f/\theta:M\to\C^n\setminus \{0\}
\]
is holomorphic and assume values in the punctured null quadric $\nullq_*=\nullq\setminus \{0\}$. 

The most important point in the development of the theory of minimal surfaces in Euclidean spaces, as presented in 
the monograph \cite{AlarconForstnericLopez2021}, is the fact that $\nullq_*$ is a complex homogeneous manifold
for the complex orthogonal group $\mathrm O_n(\C)$, hence an Oka manifold; see 
\cite[Chapter 5]{Forstneric2017E} for the latter. Thus, there is an abundance 
of holomorphic maps $h:M\to \nullq_*$ from any open Riemann surface and,  
more generally, from any Stein manifold.
Together with tools from convex integration theory and the method of dominating sprays, 
one can control the periods of the $1$-form $h\theta$ over closed curves in a given open Riemann surface $M$. 
This yields many holomorphic maps $h:M\to \nullq_*$ which integrate to conformal minimal immersions
$f:M\to\R^n$ by the Enneper--Weierstrass formula 
\[
	f(p)=f(p_0)+\int_{p_0}^p 2\Re(h\theta) \quad \text{for}\ \ p\in M,
\]
where $p_0\in M$ is a fixed reference point. Note that the integral is well-defined if and only if the 1-form
$\Re(h\theta)$ has vanishing periods. See \cite[Theorem 2.3.4]{AlarconForstnericLopez2021} for further details.

Similarly, a holomorphic immersion $f:M\to\C^n$ for $n\ge 3$ is a {\em holomorphic null curve}
if and only if 
$df=h\theta$ where $h:M\to\nullq_*$ is a holomorphic map. If $M$ is simply connected then every
conformal minimal immersion $M\to\R^n$ is the real part of a holomorphic null curve $M\to\C^n$, 
and vice versa. We can recover $f$ from $h\theta$ by the formula 
\[
	f(p)=f(p_0)+\int_{p_0}^p h\theta\quad \text{for}\ p\in M, 
\]
subject to the condition that the holomorphic 1-form $h\theta$ has vanishing periods.

%
%
\begin{definition}[Definition 1.12.9 in \cite{AlarconForstnericLopez2021}]
\label{def:admissible}
Let $M$ be a smooth surface. An {\em admissible set}\index{admissible set}
in $M$ is a compact set of the form $S=K\cup E$, where $K\subsetneq M$ is a 
finite union of pairwise disjoint compact domains in $M$ with piecewise $\Cscr^1$ boundaries
and $E = \overline{S \setminus K}$ is a union of finitely many pairwise disjoint
smooth Jordan arcs and closed Jordan curves meeting $K$ only at their endpoints
(if at all) and such that their intersections with $bK$ are transverse.
\end{definition}

%
%

We recall the notion of a generalized conformal minimal immersion and generalized holomorphic null curve.
Denote by $\Ascr^r(S,\C^n)$ the space of maps $S\to \C^n$ of class $\Cscr^r$ 
which are holomorphic in the interior of an admissible set $S$ in a Riemann surface.

\begin{definition}[Definition 3.1.2 in  \cite{AlarconForstnericLopez2021}] \label{def:GCMI}
Let $S=K\cup E$ be an admissible set in a Riemann surface $M$ (see Definition \ref{def:admissible}),
and let $\theta$ be a nowhere vanishing holomorphic $1$-form on a neighbourhood of $S$ in $M$.
A {\em generalized conformal minimal immersion} $S\to\R^n$ $(n\ge 3)$ of class $\Cscr^r$ $(r\in\N)$ 
is a pair $(f,h\theta)$, where $f: S\to \R^n$ is a $\Cscr^r$ map whose restriction to the interior $\mathring S=\mathring K$ 
is a conformal minimal immersion, and the map $h\in \Ascr^{r-1}(S,\nullq_*)$ satisfies the following two conditions:
\begin{enumerate}[\rm (a)]
\item $h\theta =2\di f$ holds on $K$, and
\item for every smooth path $\alpha$ in $M$ parameterizing a connected component of
$E=\overline{S\setminus K}$ we have that $\Re(\alpha^*(h\theta))=\alpha^*(df)=d(f\circ \alpha)$.
\end{enumerate}
\end{definition}

Note that the complex $1$-form $h\theta$ in the above definition determines 
the $1$-jet along $S$ of a conformal harmonic extension of $f$. 
With an abuse of language, we shall sometimes call the map $f$ itself a 
generalized conformal minimal immersion. The following is an analogue of this notion
for holomorphic null curves.

%
%

\begin{definition}[Definition 3.1.3 in \cite{AlarconForstnericLopez2021}] \label{def:GNC}
Let $S=K\cup E$ and $\theta$ be as in Definition \ref{def:GCMI}.
A {\em generalized holomorphic null curve} $f:S\to\C^n$ $(n\ge 3)$ of class $\Cscr^r$ $(r\in\N)$ is a pair $(f,h\theta)$
where $f\in \Ascr^r(S,\C^n)$, $h\in \Ascr^{r-1}(S,\nullq_*)$, and  the following conditions hold:
\begin{enumerate}[\rm (a)]
\item $h\theta =df=\di f$ holds on $K$ (hence $f:\mathring K\to \C^n$ is a holomorphic null curve), and
\item for any smooth path $\alpha$ in $M$ parameterizing a connected component of
$E$ we have that $\alpha^*(h\theta)=\alpha^*(df)=d(f\circ \alpha)$.
\end{enumerate}
\end{definition}

The next result 
follows immediately from \cite[Lemma 3.5.4]{AlarconForstnericLopez2021} 
and \cite[proof of Theorem 3.6.1]{AlarconForstnericLopez2021};
see the equations (3.38) and (3.39) in \cite{AlarconForstnericLopez2021}. 
The same argument applies to null curves in  part (b). The last part of the proposition is very important
in our construction.

%
%
\begin{proposition}\label{prop:extend}
Let $S=K\cup E$ be an admissible set in a Riemann surface $M$,
and let $\theta$ be a nowhere vanishing holomorphic $1$-form on 
a neighbourhood of $S$. Then, the following assertions hold for every pair of integers
$n\ge 3$ and $r\ge 1$. 
\begin{enumerate}[\rm (a)]
\item 
Every conformal minimal immersion $f:K\to\R^n$ of class $\Cscr^r$ extends 
to a generalized conformal minimal immersion $(f,h\theta)$ of class $\Cscr^r$ on $S$.
\item
Every $\Cscr^r$ map $f:K\to\C^n$ such that $f:\mathring K\to\C^n$ is a holomorphic null curve 
extends to a generalized holomorphic null curve $(f,h\theta)$ of class $\Cscr^r$ on $S$.
\end{enumerate}
If $f(K)$ is contained in a connected open set $\Omega$ in $\R^n$ or $\C^n$, respectively, 
then the extension can be chosen such that $f(S)\subset \Omega$.
More precisely, if $E$ is an arc in $M$ with the endpoints $E\cap K = \{p,q\}\in bK$ and $L$ is an arc 
in $\R^n$ or $\C^n$ connecting the points $f(p)$ and $f(q)$, then the extension of $f$
from $K$ to $K\cup E$ can be chosen such that $f(E)$ is contained in any given neighbourhood of $L$.
\end{proposition}

Recall that a compact set $S$ in an open Riemann surface $M$ is said to be {\em Runge} in $M$ if
every holomorphic function on a neighbourhood of $S$ can be approximated uniformly on $S$ 
by holomorphic functions on $M$. By Runge's theorem for Riemann surfaces,
this holds if and only if $M\setminus K$ has no relatively
compact connected components \cite[Theorem 4]{FornaessForstnericWold2020}.
The following is a simplified version of the Mergelyan approximation theorem for conformal minimal surfaces
and holomorphic null curves, given by \cite[Theorems 3.6.1 and 3.6.2]{AlarconForstnericLopez2021}. 
 
%
%
\begin{theorem}
\label{th:MergelyanCMI}
Assume that $M$ is an open Riemann surface, $S$ is an admissible Runge set in $M$, 
$n\ge 3$ and $r\ge 1$ are integers, and $f:S\to\R^n$ is a generalized conformal minimal immersion 
of class $\Cscr^r(S)$. Given $\epsilon>0$, there is a conformal minimal immersion $\tilde f:M\to \R^n$ 
satisfying $\|\tilde f-f\|_{\Cscr^r(S)}< \epsilon$. If $n\ge 5$ then $\tilde f$ can be chosen to be an injective immersion,
and if $n=4$ then $\tilde f$ can be chosen to be an immersion with simple double points.

Likewise, a generalized holomorphic null curve $f:S\to\C^n$ of class $\Cscr^r(S,\R^n)$ with $r\ge 1$
can be approximated in $\Cscr^r(S)$ by holomorphic null embeddings $M\to\C^n$.
\end{theorem}

\begin{remark}\label{rem:Mergelyan}
Since an admissible set $S$ is Runge in an open neighbourhood of itself in the ambient Riemann surface, 
Theorem \ref{th:MergelyanCMI} gives approximation of a generalized conformal minimal immersion 
on $S$ by a conformal minimal immersion on a neighbourhood of $S$.
\end{remark}

Given a compact, connected, smooth bordered surface $M$, an 
immersion $f:M\to \R^n$, and a point $p\in \mathring M$, we denote by 
$ 
	\dist_f(p,bM)
$ 
the infimum of the lengths of piecewise $\Cscr^1$ paths in $M$ connecting $p$ to $bM$
in the Riemannian metric $f^*ds^2$. This number is called the {\em intrinsic radius} of $M$
with respect to $f$. 

Our last main tool is the following lemma which says that the intrinsic radius a compact bordered Riemann surface 
can be arbitrarily large with respect to a conformal minimal immersion (or a holomorphic null curve)
which is arbitrarily uniformly close to a given one. This lemma 
is at the core of the new construction methods in the Calabi--Yau problem for minimal surfaces,
presented in \cite[Chapter 7]{AlarconForstnericLopez2021}.
The main ingredient in its proof is the Riemann-Hilbert problem for 
conformal minimal surfaces, which was first introduced in this subject in the paper
\cite{AlarconForstneric2015MA} and was developed further in 
\cite{AlarconDrinovecForstnericLopez2015PLMS,AlarconForstnericLopezMAMS}
and in \cite[Chapter 6]{AlarconForstnericLopez2021}.

%
%
\begin{lemma}\label{lem:EnlargingDiameter}
Assume that $M$ is a compact bordered Riemann surface with piecewise smooth boundary
and $f:M\to\R^n$ for $n\ge 3$ is a conformal minimal immersion of class $\Cscr^1(M)$.  
Given a point $p_0\in \mathring M$ and numbers $\epsilon>0$ (small) and $\tau>0$ (big), there is a 
conformal minimal immersion $\tilde f :M\to \R^n$ of class $\Cscr^1(M)$ such that
\[
	|\tilde f(p)-f(p)|<\epsilon\ \ \text{for all $p\in M$}\ \ \text{and}\ \ \dist_{\tilde f}(p_0,bM)>\tau.
\]
The analogous result holds for holomorphic null immersions $M\to\C^n$.
\end{lemma}

\begin{remark}\label{rem:enlarging}
Lemma \ref{lem:EnlargingDiameter} is a simplified version of 
\cite[Lemma 4.1]{AlarconDrinovecForstnericLopez2015PLMS};   
a more precise version with interpolation is \cite[Lemma 7.3.1]{AlarconForstnericLopez2021}. 
Although the boundary $bM$ is assumed to be smooth in both results, piecewise smooth boundary of
class $\Cscr^{k,\alpha}$ for some $k\in \N$ and $0<\alpha < 1$ suffices for the arguments
(see \cite[Remark 7.4.2]{AlarconForstnericLopez2021}). The intrinsic boundary distance at the 
finitely many corner points of $bM$ can be enlarged by the method of exposing points
(see \cite[Theorem 6.7.1]{AlarconForstnericLopez2021}), which is an integral part of 
the proof of \cite[Lemma 7.3.1]{AlarconForstnericLopez2021}.
\end{remark}

%
%
%
%
%
\section{Proof of Theorem \ref{th:main}}\label{sec:proof1}

We begin by recalling the standard construction of a Cantor set in the plane $\C=\R^2$.

Let $P=P_0\subset \C$ be a closed rectangle. In the first step we remove from $P$ an open horizontal 
strip of positive width around the straight line segment dividing $P$ in two halves, top and bottom.
This gives two smaller disjoint rectangles $Q_1$ and $Q_2$. Next, we remove from each of them 
an open vertical strip around the straight line segments dividing $Q_1$ and $Q_2$, respectively, in two halves, 
left and right. This gives four pairwise disjoint rectangles $P_1^j$ $(j=1,\ldots,4)$ of the same size, 
and we set $P_1=\bigcup_{j=1}^4 P_1^j$. 
Thus, the passage from $P=P_0$ to $P_1$ amounts to removing a central cross from $P$.

We now repeat the same procedure for each of the rectangles $P_1^j$, removing a central cross
in order to obtain four smaller pairwise disjoint rectangles. This gives rectangles $P_2^j$ for $j=1,2,\ldots,16$
of the second generation, and we set $P_2=\bigcup_{j=1}^{16} P_2^j$. 
Continuing inductively, we find a decreasing sequence of compacts
\begin{equation}\label{eq:Pi}
	P=P_0\supset P_1\supset P_2\supset \cdots \supset \bigcap_{i=0}^\infty P_i = C
\end{equation}
whose intersection $C$ is a Cantor set. The set $P_i$ is the union of $4^i$ pairwise disjoint closed rectangles, 
obtained by removing a central cross from each of the rectangles in $P_{i-1}$.

\begin{figure}[ht]

\begin{center}
\rotatebox{270}{\includegraphics[scale=.28]{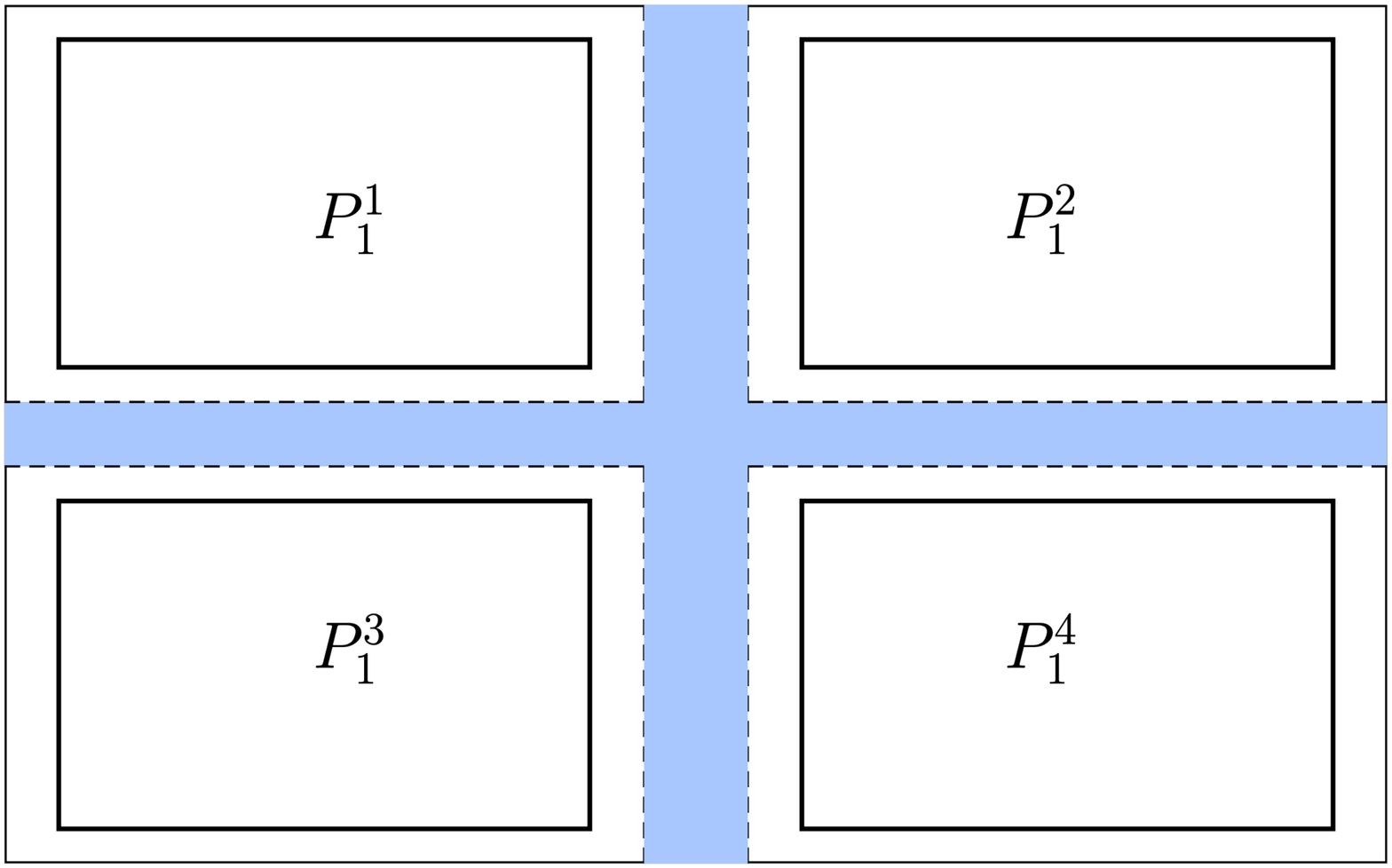}}
\end{center}

\caption{A central cross removed from a rectangle $P$}
\label{fig:slika1}

\end{figure}

In our proof of Theorem \ref{th:main}, the crosses removed at every step must be chosen fairly narrow. 
Furthermore, after removing the central cross from a given rectangle,
we will slightly shrink each of the new rectangles towards its centre in order to 
ensure that the sequence of compacts in \eqref{eq:Pi} is such that $P_{i+1}$ is 
contained in the interior of $P_i$ for every $i=0,1,2,\ldots$. 
By choosing the width of the crosses and the amount of
shrinking small enough at every step, we obtain a Cantor set $C$ whose area is 
arbitrarily close to the area of the initial rectangle $P$. The first generation of this process is shown on Figure \ref{fig:slika1}.

%
%
%
%

Theorem \ref{th:main} is a special case of the following result which we shall now prove.

%
%
\begin{theorem}\label{th:mainbis}
Assume that $M$ is a compact connected Riemann surface,
$P$ is a compact rectangle in a holomorphic coordinate chart 
on $M$, and $f:M\setminus \mathring P\to\R^n$ for $n\ge 3$ is a conformal minimal immersion 
of class $\Cscr^1(M \setminus \mathring P)$. Given an open set $\Omega \subset \R^n$ containing 
$f(M\setminus \mathring P)$, there are a Cantor set $C\subset P$ and a complete conformal 
minimal immersion $\tilde f:M\setminus C\to \Omega$ (embedding if $n\ge 5$) 
approximating $f$ as closely as desired in $\Cscr^1(M\setminus \mathring P)$. 

The analogous conclusion holds if $M=R\setminus \bigcup_i D_i$ is an open Riemann surface 
of the form \eqref{eq:bordered}, $P\subset M$ is a rectangle as above, and 
$f:R\setminus (\cup_{j=1}^m\mathring D_j \cup \mathring P)\to\R^n$ for $n\ge 3$ 
and $m\in\N$ is a conformal minimal immersion of class $\Cscr^1$ taking values in 
an open set $\Omega \subset \R^n$. In this case there exist a Cantor set $C\subset P$ 
and a continuous map $\tilde f: \overline M\setminus C\to \Omega$ such that 
$\tilde f:M\setminus C\to \Omega\subset\R^n$ is a complete conformal minimal immersion 
(embedding if $n\ge 5$) which approximates $f$ as closely as desired uniformly on 
$\overline M\setminus \mathring P$ and such that $\tilde f(bM)$ 
is a union of pairwise disjoint Jordan curves $\tilde f(bD_i)$. 
\end{theorem}

We emphasize that the Cantor set in our construction cannot be specified in advance.

\begin{proof}
For simplicity of exposition we assume that $M$ is a compact connected Riemann surface
without boundary. The case when $M$ is a bordered Riemann surface of the form \eqref{eq:bordered}  
follows by combining the procedure explained below with the one in \cite{AlarconForstneric2021RMI}.

Write $P=P_0$. 
Note that $M_0:=M\setminus \mathring P_0$ is a compact domain in $M$ with piecewise smooth boundary $bM_0=bP_0$.
Fix a base point $p_0\in \mathring M_0$ which will be used to measure the intrinsic boundary distances. 
Since $f(M\setminus \mathring P_0)$ is connected, we may assume that the neighbourhood
$\Omega$ of $f(M\setminus \mathring P_0)$ in the statement of the theorem is connected as well.

We now explain how to find the next compact domain $M_1\subset M$ with piecewise smooth boundary
satisfying $M_0 \subset \mathring M_1$, and a conformal minimal immersion $f_1:M_1\to\R^n$ which approximates 
$f_0=f$ in $\Cscr^1(M_0,\R^n)$ such that, for a given constant $c_1>0$, we have that 
\begin{equation}\label{eq:distf1}
	f_1(M_1)\subset \Omega \ \ \text{and} \ \ \dist_{f_1}(p_0,bM_1) > c_1.
\end{equation}

Let $E$ denote the horizontal straight line segment dividing the rectangle $P_0$ in top and bottom halves. 
Then, $S=M_0\cup E$ is an admissible set in $M$ (see Definition \ref{def:admissible}). 

By Proposition \ref{prop:extend} 
we can extend $f$ from $M_0$ across $E$ to a generalized conformal minimal immersion 
$f:S\to\R^n$ of class $\Cscr^1$ with $f(S)\subset\Omega$. 

By the Mergelyan approximation theorem for conformal minimal immersions 
(see Theorem \ref{th:MergelyanCMI} and Remark \ref{rem:Mergelyan}), 
we can approximate $f$ as closely as desired in $\Cscr^1(S,\R^n)$ by a conformal minimal immersion 
$g:U\to\R^n$ on an open connected neighbourhood $U\subset M$ of $S$. By shrinking $U$ around $S$ 
if necessary, we may assume that $g(U)\subset \Omega$. 

We now choose closed top and bottom rectangles $Q_1,Q_2\subset \mathring P_0 \setminus E$ 
as described above, first removing from $P_0$ a narrow horizontal strip around $E$ and then 
shrinking each of the two rectangles by a small amount, 
so that $P_0\setminus (\mathring Q_1\cup \mathring Q_2) \subset U$. Hence, the compact domain
$M'_1= M\setminus (\mathring Q_1\cup \mathring Q_2)$ in $M$ with piecewise smooth boundary 
lies in the domain $U$ of the map $g$ and it contains $M_0=M\setminus\mathring P_0$ in its interior.

Next, we repeat the same procedure with each of the two rectangles $Q_1$ and $Q_2$, 
splitting them in left and right parts by removing a narrow vertical band around their central segments $E_1$ 
and $E_2$, respectively. Using Proposition \ref{prop:extend} we extend $g$ across the arcs $E_1$ and $E_2$ 
to a generalized conformal minimal immersion on the admissible set $S'=M'_1\cup E_1\cup E_2$, 
with range in $\Omega$. By \cite[Proposition 3.3.2 (a)]{AlarconForstnericLopez2021} 
we approximate $g$ as closely as desired in $\Cscr^1(S',\R^n)$ by a conformal minimal immersion 
$\tilde g:V\to \Omega\subset \R^n$ on a neighbourhood $V$ of $S'$.
Shrinking $V$ around $S'$ we may assume that $\tilde g(V) \subset \Omega$. 
Finally, pick closed rectangles $P_1^j$ for $j=1,\ldots,4$ such that, setting 
\[
	P_1=\bigcup_{j=1}^4 P_1^j \quad\text{and}\quad M_1=M\setminus \mathring P_1,
\]
we have that $M_0\subset \mathring M_1$ and $M_1$ is contained in $V$ (the domain of $\tilde g$).

By Lemma \ref{lem:EnlargingDiameter} we can approximate  
$\tilde g$ as closely as desired in $\Cscr^1(M_0,\R^n)$ by a conformal minimal immersion 
$f_1:M_1\to \Omega\subset\R^n$ satisfying condition \eqref{eq:distf1}.

This concludes the initial step. All subsequent steps are of the same kind.
They follow the inductive construction of a Cantor set $C$, explained at the beginning of the section,
the only difference being a small shrinking of rectangles in the next generation 
obtained by removing central crosses from the rectangles of the given generation.
This shrinking was explained above when describing the initial step of the induction.

We now conclude the proof.
Pick sequences $c_i>0$ and $\epsilon_i>0$ such that $\lim_{i\to\infty}c_i=+\infty$ and
$\lim_{i\to\infty}\epsilon_i=0$. By using the above procedure we inductively construct 
\begin{itemize}
\item  a decreasing sequence of compact sets \eqref{eq:Pi} such that $P_{i+1}$ is contained in 
the interior of $P_i$ for every $i=0,1,2,\ldots$ and $C=\bigcap_{i=0}^\infty P_i$ is a Cantor set, and
\item  a sequence of conformal minimal immersions $f_i:M_i=M\setminus P_i\to \Omega\subset \R^n$, 
\end{itemize}
such that the following conditions hold for every $i=0,1,2,\ldots$: 
\begin{equation}\label{eq:ab}
	\text{(a) \ \ $\|f_{i+1}-f_i\|_{\Cscr^1(M_i)}<\epsilon_i$ \ \ \ \ and\ \ \ \ (b)\ \ $\dist_{f_i}(p_0,bM_i) > c_i$}. 
\end{equation}
Note that $\overline M_i\subset M_{i+1}$ for every $i$ and $\bigcup_{i=0}^\infty M_i = M\setminus C$.
Assuming that the numbers $\epsilon_i$ converge to zero fast enough, condition (a) ensures that the sequence
$f_i$ converges to a conformal minimal immersion 
\[
	\tilde f=\lim_{i\to\infty} f_i:M\setminus C\to\R^n
\]
whose image is contained in the given neighbourhood $\Omega$ of $f(M\setminus P_0)$.
Furthermore, if $n\ge 5$ then we may ensure that each $f_i$ and their limit $\tilde f$ are injective.
This is standard, see for example \cite[proof of Theorem 3.6.1]{AlarconForstnericLopez2021}.

Conditions (a) and (b) in \eqref{eq:ab} together clearly imply that 
\[	
	\dist_{\tilde f}(p_0,bM_i) > c_i/2\ \ \ \text{for every $i\in \N$} 
\]
provided that the sequence $\epsilon_i$ goes to zero fast enough. 
Since the increasing sequence of domains $M_i$ forms a normal exhaustion of  $M\setminus C$,
every divergent path in $M\setminus C$ emanating from $p_0$ must cross $bM_i$ 
for every $i$, and hence it has infinite length with respect to the metric $\tilde f^* ds^2$ in $M\setminus C$. 
In other words, the immersion $\tilde f$ is complete.
\end{proof}

%
%
%
%
\section{Generalization to other geometries} \label{sec:proof2}

As mentioned in the introduction, the construction in the proof of Theorem \ref{th:mainbis} 
generalizes to several other geometries listed in the following theorem.

%
%
\begin{theorem}\label{th:main2}
The analogue of Theorem \ref{th:mainbis} holds for the following classes of maps:
\begin{enumerate}[\rm (a)] 
\item Conformal harmonic immersions of nonorientable conformal surfaces into $\R^n$ for $n\ge 3$.
\item Holomorphic immersions into an arbitrary complex manifold of dimension $\ge 2$. 
\item Holomorphic null immersions into $\C^n$ with $n\ge 3$. 
\item Holomorphic Legendrian immersions into any complex contact manifold.
\item Immersed oriented superminimal surfaces in any self-dual or anti-self-dual Einstein four-manifold.
\end{enumerate}
\end{theorem}

\begin{remark}
The statement of the theorem refers to completeness of immersed surfaces with respect to a
Riemannian metric on the ambient manifold $Y$. Since the said surfaces are contained in a relatively 
compact subset of $Y$ (indeed, in a small neighbourhood of the image of the original given
compact surface), completeness does not depend on the choice of the metric on $Y$ since any two 
metrics are comparable on a relatively compact domain.
\end{remark}

\begin{proof} 
The three crucial ingredients in the proof of Theorem \ref{th:mainbis} are the following:
\begin{enumerate}[\rm (i)]
\item existence of a suitable extension of an immersion in the given class from $K$ to $S=K\cup E$, where
$S$ is an admissible set (cf.\  Proposition \ref{prop:extend}),
\item  the Mergelyan approximation theorem on admissible sets for immersions 
in the given class (see Theorem \ref{th:MergelyanCMI} for conformal minimal surfaces), and 
\item  the lemma on increasing the intrinsic radius (see Lemma \ref{lem:EnlargingDiameter}). 
\end{enumerate}
These tools are available in all geometries listed in the theorem. Let us go case by case.

\smallskip\noindent{\em Case (a): nonorientable minimal surfaces.} 
The existence of an extension from $K$ to $S=K\cup E$ is seen as in the orientable case,
the Mergelyan approximation theorem is given by 
\cite[Theorem 4.4]{AlarconForstnericLopezMAMS}, and the analogue 
of Lemma \ref{lem:EnlargingDiameter} is given by \cite[Theorem 6.6]{AlarconForstnericLopezMAMS}.

\smallskip\noindent{\em Case (b): holomorphic immersions.} 
A holomorphic immersion $f:K\to Y$ clearly extends to a smooth immersion $f:S=K\cup E\to Y$ 
with range in a given neighbourhood of $f(K)$.
The Mergelyan approximation theorem for manifold-valued maps on admissible sets
holds by \cite[Corollary 9, p.\ 178]{FornaessForstnericWold2020} and 
\cite[Theorem 1.13.1 (b)]{AlarconForstnericLopez2021}. The lemma on increasing 
the intrinsic radius was shown in \cite{AlarconForstneric2013MA} for immersions into $\C^n$ with $n\ge 2$,
which complete the proof for the case $Y=\C^n$ with $n\ge 2$. 
Combining these methods with the gluing techniques for holomorphic maps explained 
in \cite{DrinovecForstneric2007DMJ} and \cite[Chapter 5]{Forstneric2017E} implies the same 
result for immersions into any complex manifold of complex dimension at least two.

\smallskip\noindent{\em Case (c): holomorphic null immersions.} 
This case was developed in \cite{AlarconForstneric2015MA,AlarconDrinovecForstnericLopez2015PLMS}
and is covered by \cite[Theorem 7.4.12]{AlarconForstnericLopez2021}.
The proof is almost the same as the proof of Theorem \ref{th:mainbis}.

\smallskip\noindent{\em Case (d): holomorphic Legendrian immersions.} 
The case of immersions to Euclidean spaces $\C^{2n+1}$ with the standard contact form
is developed in \cite[Section 6]{AlarconForstnericLopez2017CM}.
The generalization to Legendrian immersions into an arbitrary complex contact manifold follows 
by combining \cite[Theorem 1.3]{AlarconForstneric2019IMRN} with 
the Mergelyan approximation theorem in \cite{Forstneric2020Mergelyan}.

\smallskip\noindent{\em Case (e): superminimal surfaces in (anti-) self-dual Einstein four-manifolds.} 
This follows from Case (d) by using the Bryant correspondence for Penrose twistor spaces over
such manifolds; see \cite{AlarconForstnericLarusson2021GT,Forstneric2021JGA}.
Given a self-dual or anti-self-dual Einstein four-manifold $Y$ with nonzero scalar curvature, 
the total space $X$ of the Penrose twistor bundle $\pi:X\to Y$ 
is a complex contact three-manifold, and the Bryant correspondence 
(see \cite{Bryant1982JDG,Friedrich1984}) provides a bijective correspondence between 
holomorphic and antiholomorphic immersed Legendrian curves in $X$ and oriented 
immersed superminimal surfaces in $Y$. 
In case that the Einstein metric on $Y$ has vanishing scalar curvature, the natural 
horizontal holomorphic distribution on $X$ orthogonal to the fibres of $\pi:X\to Y$ 
is integrable (a holomorphic hypersurface foliation), so 
the result reduces to that in Case (b); see \cite{Forstneric2021JGA} for the details.
\end{proof}

In this connection, we point out that Section 6 in the paper \cite{AlarconForstnericLarusson2021GT} 
describes  an axiomatic approach to the Calabi--Yau property for a given class of immersions from 
compact manifolds with boundary. If the source manifolds are surfaces and we add to those axioms
the content of Proposition \ref{prop:extend} (i.e., the existence of a suitable extension of an immersion
in  the given class from $K$ to $S=K\cup E$, where $S$ is an admissible set)
and the Mergelyan approximation theorem on admissible sets for maps in the given class, then 
the analogue of Theorem \ref{th:mainbis} holds for this class of immersions.

\noindent {\bf Acknowledgement.} 
I wish to thank Antonio Alarc\'on for an observation which led to 
an improvement of the main result, Theorem \ref{th:mainbis}.


\end{document}